\documentclass{amsart}
\usepackage{amsmath,amssymb,amscd,latexsym}

\usepackage[all]{xy}
\newcommand{\A}{{\mathbb{A}}}

\newcommand{\F}{{\mathbb{F}}}

\newcommand{\Pro}{{\mathbb{P}}}

\newcommand{\Z}{{\mathbb{Z}}}
\newcommand{\Zp}{{\mathbb{Z}}/p}
\newcommand{\Zl}{{\mathbb{Z}}/\ell}

\newcommand{\hZ}{\hat{{\mathbb{Z}}}}
\newcommand{\hFab}{\hat{F}_{\mathrm{ab}}}

\newcommand{\car}{\mathrm{char}\;}
\newcommand{\et}{\mathrm{\acute{e}t}}

\newcommand{\ok}{\overline{k}}

\newcommand{\pro}{\mathrm{pro}}

\newcommand{\Spec}{\mathrm{Spec}\,}
\newcommand{\Aut}{\mathrm{Aut}}

\newcommand{\colim}{\operatorname*{colim}}
\newcommand{\holim}{\operatorname*{holim}}
\newcommand{\cosk}{\mathrm{cosk}}
\newcommand{\Et}{\mathrm{Et}\,}
\newcommand{\hEt}{\hat{\mathrm{Et}}\,}

\newcommand{\compl}{\hat{(\cdot)}}

\newcommand{\Hom}{\mathrm{Hom}}

\newcommand{\sk}{\mathrm{sk}}

\newcommand{\Sm}{\mathrm{Sm}}

\newcommand{\sPreS}{\mathrm{sPre}(\Sm/S)}
\newcommand{\Ch}{{\mathcal C}}

\newcommand{\Eh}{{\mathcal E}}
\newcommand{\hEh}{\hat{\mathcal E}}
\newcommand{\Fh}{{\mathcal F}}
\newcommand{\Gh}{{\mathcal G}}
\newcommand{\Hh}{{\mathcal H}}
\newcommand{\Hhf}{{\mathcal H}_{\mathrm{fin}}}

\newcommand{\hHh}{\hat{{\mathcal H}}}

\newcommand{\Mh}{{\mathcal M}}

\newcommand{\Rh}{{\mathcal R}}
\newcommand{\Sh}{{\mathcal S}}

\newcommand{\hSHh}{{\hat{\mathcal{SH}}}}

\newcommand{\hSh}{\hat{\mathcal S}}

\newcommand{\hShp}{\hat{\mathcal S}_{\ast}}

\newcommand{\hSp}{\mathrm{Sp}(\hShp)}

\newcommand{\Sp}{\mathrm{Sp}(\Sh_{\ast})}

\newcommand{\Th}{{\mathcal T}}
\newcommand{\Xh}{\mathcal{X}}

\newtheorem{theorem}{Theorem}[section]
\newtheorem{lemma}[theorem]{Lemma}
\newtheorem{prop}[theorem]{Proposition}
\newtheorem{defn}[theorem]{Definition}

\newtheorem{cor}[theorem]{Corollary}

\newtheorem{remark}[theorem]{Remark}
\begin{document}
\title{Profinite Homotopy Theory}
\author{Gereon Quick}
\date{}
\begin{abstract}
We construct a model structure on simplicial profinite sets such that the homotopy groups carry a natural profinite structure. This yields a rigid profinite completion functor for spaces and pro-spaces. One motivation is the \'etale homotopy theory of schemes in which higher profinite \'etale homotopy groups fit well with the \'etale fundamental group which is always profinite. We show that the profinite \'etale topological realization functor is a good object in several respects. 
\end{abstract}
\maketitle

\section{Introduction}

Let $\Sh$ be the category of simplicial sets and let $\Hh$ be its homotopy category. The \'etale homotopy type of a scheme had been defined by Artin and Mazur \cite{artinmazur} as a pro-object in the homotopy category $\Hh$. Friedlander \cite{fried} rigidified the definition such that the \'etale topological type is a pro-object in the category of simplicial sets itself. This theory has found very remarkable applications. Quillen and Friedlander used it to prove the Adams conjecture, a purely topological problem. Friedlander defined and studied \'etale K-theory, an \'etale topological analogue of algebraic K-theory. More recently, Dugger and Isaksen \cite{sumsofsquares} proved a sums-of-squares formula over fields of positive characteristic using \'etale topological arguments. Schmidt \cite{ringsofintegers} answered open questions in Galois theory using \'etale homotopy groups. In \cite{etalecob}, a stable \'etale realization functor has been constructed and an \'etale cobordism theory for schemes has been defined. For smooth schemes over an algebraically closed field and with finite coefficients, \'etale cobordism agrees with algebraic cobordism after inverting a Bott element.\\
In almost every application, one considers a profinitely completed object, either with respect to finite or even with respect to $p$-groups for some prime $p$. The profinite (resp.~pro-$p$-) completion is an object that is universal with respect to maps to spaces whose homotopy groups are finite (resp.~$p$-) groups. Artin and Mazur realized the profinite completion of a space or pro-space only in pro-$\Hh$. There is no profinite completion for Friedlander's rigid objects in pro-$\Sh$. It is well known that in many respects it is preferable to work in a model category itself and not only in the corresponding homotopy category. Hence it is a fundamental question if there is a rigid model for the profinite completion, i.e. an object in $\Sh$ which is homotopy equivalent to the Artin-Mazur completion in pro-$\Hh$. Since the \'etale fundamental group of a scheme is always a profinite group, this is equivalent to the fundamental question if there is a space, not only a pro-object in some homotopy category, that yields profinite higher \'etale homotopy groups.\\
Bousfield-Kan \cite{bouskan} proved the existence of a $\Z/p$-completion for a simplicial set for every prime $p$ in $\Sh$. Moreover, Morel \cite{ensprofin} showed that there is even a $\Z/p$-model structure on the category $\hSh$ of simplicial profinite sets such that the Bousfield-Kan completion yields a fibrant replacement functor in $\hSh$. In \cite{etalecob}, this rigid model for the pro-$p$-homotopy theory has been used.  In this paper we prove that there is a model for arbitrary profinite completion. We do this by constructing a suitable model structure on the category of simplicial profinite sets such that the homotopy groups carry a natural profinite structure.\\
The plan of this paper is the following. First, we construct profinite fundamental groups and then we use the profinite topology to define continuous cohomology with local topological coefficients for profinite spaces. The main technical result is this: There is a model structure on simplicial profinite sets such that a weak equivalence is a map that induces isomorphisms on fundamental groups and in cohomology with local finite abelian coefficients. This model structure is fibrantly generated, simplicial and left proper. This result enables us to define higher profinite homotopy groups. It is an important property of the category of profinite spaces that the limit functor is homotopy invariant for cofiltering diagrams.\\
The fibrant replacement of the completion of a simplicial set in $\hSh$ is a rigid model for the Artin-Mazur profinite completion of \cite{artinmazur}, and is equivalent in an appropriate sense to the completions of Bousfield-Kan \cite{bouskan} and Morel \cite{ensprofin}. An important advantage of this approach is that, for profinite local coefficients, the continuous cohomology of this completion coincides with the continuous cohomology of Dwyer-Friedlander \cite{dwyfried}.\\
Although we had been motivated by \'etale homotopy theory, we have postponed this application to the last section of the paper. As in \cite{etalecob}, we consider a profinite \'etale topological type of a locally noetherian scheme based on the work of Artin-Mazur and Friedlander. The resulting profinite space has three main advantages. First, its cohomology agrees with the continuous \'etale cohomology for a locally constant profinite sheaf defined by Dwyer-Friedlander in \cite{dwyfried} and by Jannsen in \cite{jannsen}. Second, it provides a rigid model for the profinite higher \'etale homotopy groups of a scheme first defined in \cite{artinmazur}. Third, it can be used for an \'etale realization functor of the motivic stable homotopy category. We discuss this last point briefly at the end of the paper, indicating that the \'etale realization also yields a derived functor from the flasque model structure of \cite{flasque}.\\
{\bf Acknowledgments}: I would like to thank Alexander Schmidt and the referee for a lot of detailed and helpful comments. 

\section{Homotopy Theory of Profinite Spaces}

For a category $\Ch$ with small limits, the pro-category of $\Ch$, denoted pro-$\Ch$, has as objects all cofiltering diagrams $X:I \to \Ch$. Its sets of morphisms are defined as
$$\Hom_{\pro-\Ch}(X,Y):=\lim_{j\in J}\colim_{i\in I} \Hom_{\Ch}(X_i,Y_j).$$
A constant pro-object is one indexed by the category with one object and one identity map. The functor sending an object $X$ of $\Ch$ to the constant pro-object with value $X$ makes $\Ch$ a full subcategory of pro-$\Ch$. The right adjoint of this embedding is the limit functor $\lim$: pro-$\Ch$ $\to \Ch$, which sends a pro-object $X$ to the limit in $\Ch$ of the diagram corresponding to $X$.\\
Let $\Eh$ denote the category of sets and let $\Fh$ be the full subcategory of finite sets. Let $\hEh$ be the category of compact and totally disconnected topological spaces. We may identify $\Fh$ with a full subcategory of $\hEh$ in the obvious way. The limit functor $\lim$: pro-$\Fh \to \hEh$ is an equivalence of categories.\\
We denote by $\hSh$ (resp. $\Sh$) the category of simplicial objects in $\hEh$ (resp. simplicial sets). The objects of $\hSh$ (resp. $\Sh$) will be called {\em profinite spaces} (resp. {\em spaces}). The forgetful functor $\hEh \to \Eh$ admits a left adjoint $\compl:\Eh \to \hEh$. It induces a functor $\compl:\Sh \to \hSh$, which is called {\em profinite completion}. It is left adjoint to the forgetful functor $|\cdot|:\hSh \to \Sh$ which sends a profinite space to its underlying simplicial set.\\ 
For a profinite space $X$ we define the set $\Rh(X)$ of simplicial open equivalence relations on $X$. An element $R$ of $\Rh(X)$ is a simplicial profinite subset of the product $X\times X$ such that, in each degree $n$, $R_n$ is an equivalence relation on $X_n$ and an open subset of $X_n\times X_n$. It is ordered by inclusion. For every element $R$ of $\Rh(X)$, the quotient $X/R$ is a simplicial finite set and the map $X \to X/R$ is a map of profinite spaces. The canonical map $X \to \lim_{R\in \Rh(X)} X/R$ is an isomorphism in $\hSh$, cf. \cite{ensprofin}, Lemme 1. Nevertheless, Isaksen pointed out that $\hSh$ is not equivalent to the category of pro-objects of finite simplicial sets.\\
Let $X$ be a profinite space. The continuous cohomology $H^{\ast}(X;\pi)$ of $X$ with coefficients in the topological abelian group $\pi$ is defined as the cohomology of the complex $C^{\ast}(X;\pi)$ of continuous cochains of $X$ with values in $\pi$, i.e. $C^n(X;\pi)$ denotes the set $\Hom_{\hat{\Eh}}(X_n,\pi)$ of continuous maps $\alpha:X_n \to \pi$ and the differentials $\delta^n:C^n(X;\pi)\to C^{n+1}(X;\pi)$ are the morphisms associating to $\alpha$ the map $\sum_{i=0}^{n+1}\alpha \circ d_i$, where $d_i$ denotes the $i$th face map of $X$. If $\pi$ is a finite abelian group and $Z$ a simplicial set, then the cohomologies $H^{\ast}(Z;\pi)$ and $H^{\ast}(\hat{Z};\pi)$ are canonically isomorphic.\\
If $G$ is an arbitrary profinite group, we may still define the first cohomology of $X$ with coefficients in $G$ as done by Morel in \cite{ensprofin} p. 355. The functor $X\mapsto \Hom_{\hEh}(X_0,G)$ is represented in $\hSh$ by a profinite space $EG$. We define the $1$-cocycles $Z^1(X;G)$ to be the set of continuous maps $f:X_1 \to G$ such that $f(d_0x)f(d_2x)=f(d_1x)$ for every $x \in X_1$. The functor $X\mapsto Z^1(X;G)$ is represented by a profinite space $BG$. Explicit constructions of $EG$ and $BG$ may be given in the standard way, see \cite{may}. Furthermore, there is a map $\delta:\Hom_{\hSh}(X,EG) \to Z^1(X;G)\cong \Hom_{\hSh}(X,BG)$ which sends $f:X_0 \to G$ to the $1$-cocycle $x\mapsto \delta f(x)=f(d_0x)f(d_1x)^{-1}$. We define $B^1(X;G)$ to be the image of $\delta$ in $Z^1(X;G)$ and we define the pointed set $H^1(X;G)$ to be the quotient $Z^1(X;G)/B^1(X;G)$. Finally, if $X$ is a profinite space, we define $\pi_0X$ to be the coequalizer in $\hEh$ of the diagram $d_0,d_1:X_1 \rightrightarrows X_0$.   

\subsection{Profinite fundamental groups} 
For the definition of a profinite fundamental group of a profinite space $X$, we follow the ideas of Grothendieck in \cite{sga1}. As in \cite{gz}, we denote by $R/X$ the category of coverings of $X$, i.e. the full subcategory of $\hSh/X$ whose objects are the morphisms $p:E\to X$ such that for each commutative diagram
\begin{equation}
\xymatrix{
\Delta[0] \ar[d]_i \ar[r]^u & E \ar[d]^p \\
\Delta[n] \ar[r]_v & X }
\end{equation}
there is a unique morphism $s:\Delta[n] \to E$ satisfying $p\circ s=v$ and $s \circ i=u$. \\
Let $g:X' \to X$ be a morphism of $\hSh/X$ and $p:E\to X$ be a covering of $X$. The cartesian square
\begin{equation}
\xymatrix{
X'\times_X E \ar[d]_{p'} \ar[r] & E \ar[d]^p \\
X' \ar[r]_g & X} 
\end{equation}
defines a covering $p'$ of $X'$. The correspondence $R/g: R/X \to R/X'$ so defined is a covariant functor.\\
Let $x\in X_0$ be a vertex of $X$ and let $\tilde{x}:\Delta[0]\to X$ be the corresponding continuous map. Since there is only one morphism $\Delta[n] \to \Delta[0]$ for each $n$, it follows from the definitions that each simplex of a covering $E$ of $\Delta[0]$ is determined by one of its vertices, and hence  $E=\sk_0E$. Thus we can identify the category $R/\Delta[0]$ with the category of profinite sets. The functor $R/\tilde{x}$ is the fiber functor of $R/X$ over $x$ and we consider its range in $\hEh$. As in \cite{gz} App.~I 2.2, it is easy to show that a morphism $p:E \to X$ in $\hSh$ is a covering if and only if it is a locally trivial morphism whose fibers are constant simplicial profinite sets.\\
In order to define the fundamental group of $X$ at $x \in X_0$ we consider the full subcategory $R_f/X$ of $R/X$ of coverings with finite fibers together with the fiber functor $R_f/\tilde{x}$. We call an object of $R_f/X$ a finite covering of $X$. The pair $(R_f/X,R/\tilde{x})$ obviously satisfies the axioms $(G1)$ to $(G6)$ of \cite{sga1} Expose V \S 4. Hence $R/\tilde{x}$ is pro-representable by a pro-object $(X_i)_{i \in I}$ in $R_f/X$, where the $X_i$ are Galois, i.e. the action of $\Aut(X_i)$ on the fiber of $X_i$ at $x$ is simply transitive. The corresponding limit of this pro-object defines an element $\tilde{X}$ in $R/X$ of profinite coverings of $X$ and may be considered as the universal covering of $X$ at $x$. We denote the automorphism group of $\tilde{X}$ by $\pi_1(X,x)$ and call it {\em the profinite fundamental group of $X$ at $x$}. It has a canonical profinite structure as the limit of the finite automorphism groups $\Aut(X_i)$.\\ 
For varying $x$ and $y$ in $X_0$, we have two corresponding fiber functors $R/x$ and $R/y$ and representing objects $(\tilde{X},x)$ and $(\tilde{X},y)$. The morphisms between them are in fact isomorphisms. Hence we may consider the fundamental groupoid $\Pi X$ of $X$ whose objects are the vertices of $X$ and whose morphisms are the sets $\Hom((\tilde{X},x),(\tilde{X},y))$. 
Now let $X$ be a pointed simplicial set and let $\hat{X}$ be its profinite completion. The finite coverings of $X$ and $\hat{X}$ agree as coverings in $\Sh$. Since the automorphism group of a finite covering of $X$ corresponds to a finite quotient of $\pi_1(X)$ and since $\pi_1(\hat{X})$ is the limit over these finite groups, we deduce the following result.
\begin{prop}\label{profinitecompletion}
For a pointed simplicial set $X$, the canonical map from the profinite group completion of $\pi_1(X)$ to $\pi_1(\hat{X})$ is an isomorphism, i.e. $\widehat{\pi_1(X)}\cong \pi_1(\hat{X})$ as profinite groups. 
\end{prop}

Moreover, we get a full description of coverings in $\hSh$. 
\begin{prop}
Let $X$ be a connected pointed profinite space with fundamental group $\pi:=\pi_1(X)$. Then the functor sending a profinite covering to its fiber is an equivalence between the category of profinite coverings of $X$ and the category of profinite sets with a continuous $\pi$-action. Its inverse is given by $S\mapsto P\times_{\pi}S$.
\end{prop}
\begin{proof}
Since $X$ is connected, a covering of $X$ is determined up to isomorphism by its fiber as a $\pi$-set. Hence the category of pro-objects of $R_f/X$ is equivalent to the category $R/X$ of profinite coverings, i.e. the limit functor $\mathrm{pro-}R_f/X \stackrel{\sim}{\to}R/X$ is an equivalence. Now the assertion follows from \cite{sga1}, Expos\'e V, Th\'eor\`eme 4.1 and Corollaire 5.9.
\end{proof}
\begin{cor}\label{classificationofcoverings}
Let $X$ be a connected pointed profinite space. There is a bijective correspondence between the sets of profinite coverings of $X$ and closed subgroups of $\pi_1(X)$.
\end{cor}
\subsection{Local coefficient systems}
Let $\Gamma$ be a profinite groupoid, i.e. a small category whose morphisms are all isomorphisms and whose set of objects and morphisms are profinite sets. A natural example is $\Pi X$ for a profinite space $X$. A (topological) local coefficient system $\Mh$ on $\Gamma$ is a contravariant functor from $\Gamma$ to topological abelian groups. A morphism $(\Mh,\Gamma) \to (\Mh',\Gamma')$ of local coefficient systems is a functor $f:\Gamma \to \Gamma'$ and natural transformation $\Mh' \to \Mh\circ f$. We call $\Mh$ a profinite local coefficient system if it takes values in the category of profinite abelian groups. A local coefficient system $\Mh$ on a profinite space $X$ is a coefficient system on the fundamental groupoid $\Gamma=\Pi X$ of $X$.\\
In order to define the continuous cohomology of a profinite space with local coefficients we use one alternative characterization of cohomology with local coefficients of a simplicial set following Goerss and Jardine \cite{gj}. We recall briefly the constructions of \cite{gj} VI \S 4. Let $\Gamma$ be a profinite groupoid and let $B\Gamma$ be its profinite classifying space. If $\gamma$ is an object of $\Gamma$, we form the category $\Gamma \downarrow \gamma$ whose objects are morphisms $\gamma' \to \gamma$ and whose morphisms are commutative diagrams in $\Gamma$. The forgetful functor $\Gamma\downarrow \gamma \to \Gamma$ induces a map $\pi_y:B(\Gamma\downarrow \gamma) \to B\Gamma$. A map $\gamma_1 \to \gamma_2$ induces a functor $B(\Gamma\downarrow \gamma_1) \to B(\Gamma\downarrow \gamma_2)$ which commutes with the forgetful functor. Hence we get a functor $\Gamma \to \hSh$ defined by $\gamma \mapsto B(\Gamma\downarrow \gamma)$ which is fibered over $B\Gamma$.\\
Now let $\Phi:X \to B\Gamma$ be a map in $\hSh$. One defines a collection of spaces $\tilde{X}_{\gamma}$ by forming pullbacks
$$\xymatrix{
\tilde {X}_{\gamma} \ar[d] \ar[r] & B(\Gamma\downarrow \gamma) \ar[d] \\
X \ar[r] & B\Gamma}$$
and gets a functor $\Gamma \to \hSh$, called covering system for $\Phi$. If $\Gamma$ is $\Pi X$ and $\Phi$ is the canonical map, then $\tilde{X}_{\gamma}$ is just the universal covering $(\tilde{X},x)$ for $x=\gamma$.\\
Moreover, if $Y:\Gamma \to \hSh$ is a functor and $\Mh$ is a topological local coefficient system, there is a corresponding cochain complex $\hom_{\Gamma}(Y,\Mh)$, having $n$-cochains given by the group of $\hom_{\Gamma}(Y_n,\Mh)$ of all {\em continuous} natural transformations, i.e. for every $\gamma \in \Gamma$ we consider only the continuous maps $Y_n(\gamma) \to \Mh(\gamma)$ functorial in $\gamma$. The differentials are given by the alternating sum of the face maps as above. For $Y=\tilde{X}$, we denote this cochain complex by $C_{\Gamma}^{\ast}(X,\Mh):=\hom_{\Gamma}(\tilde{X},\Mh)$.
\begin{defn}
For a topological local coefficient system $\Mh$ on $\Gamma$ and a map $X\to B\Gamma$ in $\hSh$, we define the {\em continuous cohomology of $X$ with coefficients in $\Mh$}, denoted by $H_{\Gamma}^{\ast}(X,\Mh)$, to be the cohomology of the cochain complex $C_{\Gamma}^{\ast}(X,\Mh)$.
\end{defn}
If $\Mh$ is a local system on $\Gamma=\Pi X$, we write $H^{\ast}(X,\Mh)$ for $H_{\Pi X}^{\ast}(X,\Mh)$. If $G$ is a profinite group and $X=BG$ is its classifying space, then a topological local coefficient system $\Mh$ on $BG$ corresponds to a topological $G$-module $M$. The continuous cohomology of $BG$ with coefficients in $\Mh$ equals then the continuous cohomology of $G$ with coefficients in $M$ as defined in \cite{tate}, i.e. 
$$H^{\ast}(BG, \Mh)=H^{\ast}(G,M).$$
The following proposition is a standard result. We restate it in our setting of profinite groups and continuous cohomology. The proof follows the arguments for \cite{mcc} Theorem $8^{\mathrm{bis}}.9$.
\begin{prop}\label{cartanlerayss}
Let $\Gamma$ be a connected profinite groupoid, let $X \to B\Gamma$ be a profinite space over $B\Gamma$ and let $\Mh$ be a local coefficient system on $\Gamma$ of topological abelian groups. For an object $\gamma$ of $\Gamma$ we denote by $\pi$ the profinite group $\Hom_{\Gamma}(\gamma,\gamma)$. Let $p:\tilde{X} \to X$ be the corresponding covering space. Then $\pi$ acts continuously on the discrete abelian group $H^q(\tilde{X};p^{\ast}\Mh)$ and there is a strongly converging Cartan-Leray spectral sequence
$$E_2^{p,q}=H^p(\pi;H^q(\tilde{X};p^{\ast}\Mh)) \Rightarrow H^{p+q}_{\Gamma}(X;\Mh).$$
\end{prop}
\begin{proof}
Since $\Gamma$ is connected, $\Mh$ corresponds to a profinite $\pi$-module $M$. Moreover, the set $C^q(\tilde{X};p^{\ast}\Mh)$ equals the set of continuous cochains $C^q(\tilde{X};M)$ denoted by $C^q$ for short. We denote by $F^{\ast} \to \hZ \to 0$ a free profinite resolution of $\hZ$ by right $\pi$-modules, where $\hZ$ is considered as a trivial $\pi$-module. We define the double complex $\Ch^{\ast,\ast}$ by
$\Ch^{p,q}:=\Hom(F^p, C^q), ~ \delta_F \otimes 1 + (-1)^q1\otimes \delta_{C}$, where the Hom-set is taken in the category of continuous $\pi$-modules.\\ 
First, we filter $\Ch^{\ast,\ast}$ row-wise and get for $E_0^{\ast,q}=\Hom_{\pi-\mathrm{modules}}(F^{\ast}, C^q)$, the complex computing $H^{\ast}(\pi;C^q)$. Since $C^q$ is a free $\pi$-module, the $E_1$-terms are concentrated in the $0$-column, where we get $H^0(\pi,C^q)=C^q(\tilde{X};M)^{\pi}$, the $\pi$-fixed points of $C^q$. By \cite{mcc} Theorem 2.15, the $E_2$-terms for this filtration are the continuous $\pi$-equivariant cohomology groups $H^{\ast}_{\pi}(\tilde{X};M)$ and the spectral sequence degenerates at $E_2$. Since $\Gamma$ is connected, these cohomology groups are canonically isomorphic to $H^{p+q}_{\Gamma}(X;\Mh)$.\\
Now we filter $\Ch^{\ast,\ast}$ column-wise and we get $E_0^{p,\ast}=\Hom_{\pi}(F^p, C^{\ast})$. Since $F^p$ is a free $\pi$-module, we may identify the corresponding $E_1^{p,\ast}$ with the complex computing $H^{\ast}_{\mathrm{cont}}(\pi;H^{\ast}(\tilde{X};M))$. Both spectral sequences strongly converge to the same target and the assertion is proved. 
\end{proof}

\subsection{The model structure on profinite spaces}
\begin{defn}\label{defnwe}
A morphism $f:X\to Y$ in $\hSh$ is called a {\em weak equivalence} if the induced map $f_{\ast}:\pi_0(X) \to \pi_0(Y)$ is an isomorphism of profinite sets, for every vertex $x\in X_0$ the map 
$f_{\ast}:\pi_1(X,x) \to \pi_1(Y,f(x))$ is an isomorphism of profinite groups and $f^{\ast}:H^q(Y;\Mh) \to H^q(X;f^{\ast}\Mh)$ is an isomorphism for every local coefficient system $\Mh$ of finite abelian groups on $Y$ for every $q\geq 0$.
\end{defn}
We want to show that this class of weak equivalences fits into a simplicial fibrantly generated left proper model structure on $\hSh$. For every natural number $n\geq 0$ we choose a finite set with $n$ elements, e.g. the set $\{0,1, \ldots, n-1\}$, as a representative of the isomorphism class of sets with $n$ elements. We denote the set of these representatives by $\Th$. Moreover, for every isomorphism class of finite groups, we choose a representative with underlying set $\{0,1,\ldots,n-1\}$. Hence for each $n$ we have chosen as many groups as there are relations on the set $\{0,1,\ldots,n-1\}$. This ensures that the collection of these representatives forms a set which we denote by $\Gh$.\\
Let $P$ and $Q$ be the following two sets of morphisms:
$$\begin{array}{lll}
P & \mathrm{consisting~of} & EG \to BG, BG \to \ast,~ L(M,n)\to K(M,n+1),\\
   &   &  K(M,n)\to \ast,~K(S,0) \to \ast\\
    &  & \mathrm{for~every~finite~set}~S\in \Th,~\mathrm{every~finite~group}~G\in \Gh,\\
   &   &  \mathrm{every~finite~abelian~group}~M\in \Gh ~\mathrm{and~every}~ n\geq 0;\\
Q & \mathrm{consisting~of} & EG \to \ast,~ L(M,n) \to \ast~ \mathrm{for~every~finite~group}~G\in \Gh,\\
   &   &  \mathrm{every~finite~abelian~group}~M\in \Gh ~\mathrm{and~every}~ n\geq 0. 
\end{array}$$
\begin{lemma}\label{injectivity}
The underlying profinite set of a profinite group $G$ is an injective object in $\hEh$. 
\end{lemma}
\begin{proof}
This can be deduced from Proposition 1 of \cite{serre}. Let $X\hookrightarrow Y$ be a monomorphism in $\hSh$ and let $f:X \to G$ be a map in $\hEh$. Since finite sets are injective objects in $\hEh$ by \cite{ensprofin} Lemme 2, there is a lift $Y \to G/U$ for every open (and hence closed) normal subgroup $U$ of $G$. Let $N$ be the set of pairs $(S,s)$ of closed subgroups $S$ of $G$ such that there is a lift $s:Y\to G/S$ making the diagram
$$\xymatrix{
X \ar[d]_i  \ar[r]^{f/S} & G/S \\
Y \ar[ur]_s}$$ 
commute. The set $N$ contains the open normal subgroups and has a natural ordering given by $(S,s)\geq (S',s')$ if $S \subseteq S'$ and $s:Y\to G/S$ is the composite of $s'$ and $G/S' \to G/S$. As shown in \cite{serre}, $N$ is an inductively ordered set and has a maximal element by Zorn's Lemma. We have to show that a maximal element $(S,s)$ of $N$ satisfies $S=\{1\}$.\\
Suppose $S\neq \{1\}$. Then there is an open subgroup $U$ of $G$ such that $S\cap U \neq S$ and $S/S\cap U$ is a finite group. By \cite{serre}, Proposition 1, $G/S\cap U$ is isomorphic as a profinite set to the product $G/S \times S/S\cap U$. The map $f/S\cap U:X\to G/S\cap U\cong G/S \times S/S\cap U$ induces a compatible map $X \to S/S\cap U$. Since $S/S\cap U$ is a finite set, it is an injective object in $\hEh$ and there is a lift $t:Y\to S/S\cap U$. Hence $s$ and $t$ define a lift $\tilde{s}:Y\to G/S\cap U$ in contradiction to the maximality of $(S,s)$. Hence $S$ is trivial and there is a lift of the initial map $f$.\\ 
\end{proof}
\begin{lemma}\label{lemme2}
1) The morphisms in $Q$ have the right lifting property with respect to all monomorphisms. 2) The morphisms in $P$ have the right lifting property with respect to all monomorphisms that are also weak equivalences.
\end{lemma}
\begin{proof}
1) Let $X\hookrightarrow Y$ be a monomorphism in $\hSh$. We have to show that every $X\to EG$, resp. $X \to L(M,n)$, can be lifted to a map $Y\to EG$, resp. $Y\to L(M,n)$, in $\hSh$. Hence we must prove that every map $X_n \to G$ may be lifted to a map $Y_n \to G$ in $\hEh$ for every profinite group $G$. This follows from Lemma \ref{injectivity}.\\
2) Let $X\hookrightarrow Y$ be a monomorphism in $\hSh$ that is also a weak equivalence. Let $G$ be a finite group, which is supposed to be abelian if $n\geq2$. We know by 1) that the morphism of complexes $C^{\ast}(Y;G)\to C^{\ast}(X;G)$ is surjective. By assumption, it also induces an isomorphism on the cohomology. Hence the maps $Z^n(Y;G) \to Z^n(X;G)$ and
$C^n(Y;G) \to C^n(X;G)\times_{Z^{n+1}(X;G)}Z^{n+1}(Y;G)$ are surjective and the maps $L(G,n) \to K(G,n+1)$ and $K(G,n) \to \ast$ have the desired right lifting property.\\
For $K(S,0) \to \ast$ and a given map $X\to K(S,0)$, we recall that $K(S,0)_n$ is equal to $S$ in each dimension and all face and degeneracy maps are identities. Hence a map $X\to K(S,0)$ in $\hSh$ is completely determined by its values on $\pi_0X$. Since $f$ induces an isomorphism on $\pi_0$, there is a lift $Y \to K(S,0)$. Hence $K(S,0)\to \ast$ also has the desired right lifting property.
\end{proof}

We remind the reader of the following definitions of \cite{ensprofin} p. 360. Let $G$ be a simplicial profinite group and let $E$ be a profinite $G$-space. We say that $E$ is a principal profinite $G$-space if, for every $n$, the profinite $G_n$-set $E_n$ is free.\\
A {\em principal $G$-fibration with base $X$} is a profinite $G$-space $E$ and a morphism $f:E\to X$ that induces an isomorphism $E/G \cong X$. We denote by $\Phi^G(X)$ the set of isomorphism classes of principal $G$-fibrations with base $X$. The correspondence $X\mapsto \Phi^G(X)$ defines via pullback a contravariant functor $\hSh^{\mathrm{op}} \to \Eh$.  
\begin{lemma}
Let $X$ be a connected profinite space and let $x\in X_0$ be a vertex. Let $G$ be a profinite group and let $\Hom(\pi_1(X,x),G)_G$ be the set of outer continuous homomorphisms. Then we have a natural isomorphism 
$$\phi:H^1(X;G)\stackrel{\cong}{\longrightarrow} \Hom(\pi_1(X,x),G)_G.$$ 
\end{lemma}
\begin{proof}
Recall that $H^1(X;G)$ equals the quotient of $\Hom_{\hSh}(X,BG)$ modulo maps that are induced by the canonical principal $G$-fibration $EG \to BG$. The map $\phi$ is defined as follows: Given a map $f:X\to BG$ we consider the induced map $X\times_{G}EG \to X$. This is a covering on which $G$ acts freely. By Corollary \ref{classificationofcoverings}, there is a quotient $Q$ of $\pi=\pi_1(X,x)$ such that $Q$ acts on the fibre $G$ of $X\times_{G}EG \to X$. This action defines a homomorphism of profinite groups $\pi \to G$ up to inner automorphisms.\\
We define an inverse $\psi$ for $\phi$. Let $\alpha: \pi \to G$ be a homomorphism of profinite groups up to inner automorphisms. Again by Corollary \ref{classificationofcoverings}, there is a covering  $X(G) \to X$ on which $G$ acts freely. Since its total space is cofibrant and $EG \to \ast$ has the left lifting property with respect to all cofibrations, there is a map $X(G) \to EG$. The corresponding quotient map is the map $\psi(\alpha):X\to BG$. One can easily check that $\phi$ and $\psi$ are mutually inverse to each other. 
\end{proof}

The proof of the previous lemma also explains the following classification of principal $G$-fibrations.
\begin{prop}\label{principalfibrations}
Let $G$ be a simplicial profinite group. For any profinite space $X$, the map
$$\theta:H^1(X;G) \to \Phi^G(X),$$
sending the image of $f:X\to BG$ in $H^1(X;G)$ to the pullback of $EG \to BG$ along $f$, is a bijection.
\end{prop}

\begin{prop}\label{whitehead}
Let $f:X\to Y$ be a map of profinite spaces. If $x \in X$ is a $0$-simplex, let $y=f(x)$ and $p:(\tilde{X},\tilde{x})\to (X,x)$, resp. $q:(\tilde{Y},\tilde{y})\to (Y,y)$, be the universal coverings and $\tilde{f}:\tilde{X}\to \tilde{Y}$ be the unique covering of $f$ with $\tilde{f}(\tilde{x})=\tilde{y}$.
The following assertions are equivalent:\\
1) The map $f$ is a weak equivalence in $\hSh$ in the sense of Definition \ref{defnwe}.\\
2) The induced maps $f^0:H^0(Y;S) \to H^0(X;S)$ for every finite set $S$,\\ 
$f^1:H^1(Y;G) \to H^1(X;G)$ for every every finite group $G$ and $f^{\ast}:H^q(Y;\Mh) \to H^q(X;f^{\ast}\Mh)$ for every local coefficient system $\Mh$ of finite abelian groups on $Y$ for every $q\geq 0$ are all isomorphisms.\\
3) The induced map $f_{\ast}:\pi_0(X) \to \pi_0(Y)$ is an isomorphism of profinite sets, $f_{\ast}:\pi_1(X,x) \to \pi_1(Y,f(x))$ is an isomorphism of profinite groups and the maps $\tilde{f}^{\ast}:H^q((\tilde{Y},f(x));M) \to H^q((\tilde{X},x);M)$ are isomorphisms for every finite abelian group $M$ for every $q\geq 0$ and every $0$-simplex $x\in X_0$.
\end{prop}
\begin{proof}
From $H^0(X;S)= \Hom_{\hEh}(\pi_0(X),S)$ for every finite set $S$, we conclude that $\pi_0(f)$ is an isomorphism if and only if $H^0(f,S)$ is an isomorphism for every finite set $S$. From the previous lemma we get that $\pi_1(f)$ is an isomorphism if and only if $H^1(f,G)$ is an isomorphism for every finite group $G$. Hence (1) and (2) are equivalent.\\
In order to show $(3) \Rightarrow (1)$, we may assume that $X$ and $Y$ are connected profinite spaces and that $f_{\ast}:\pi_1(X,x) \to \pi_1(Y,f(x))$ is an isomorphism of profinite groups for every vertex $x\in X_0$. Let $x_0 \in X$ be a fixed $0$-simplex and set $\pi:=\pi_1(X,x_0)=\pi_1(Y,f(x_0))$. If $\Mh$ is a local coefficient system on $Y$, then there is a morphism of Cartan-Leray spectral sequences of Proposition \ref{cartanlerayss}
$$\begin{array}{cccc}
E_2^{p,q}= & H^p(\pi,H^q(\tilde{Y},q^{\ast}\Mh)) & \Longrightarrow & H^{p+q}(Y,\Mh)\\
  & \downarrow &  & \downarrow \\
E_2^{p,q}= & H^p(\pi,H^q(\tilde{X},p^{\ast}f^{\ast}\Mh)) & \Longrightarrow & H^{p+q}(X,f^{\ast}\Mh).
\end{array}$$
If we assume (3), then the map on the $E_2$-terms is an isomorphism and we get an isomorphism on the abutments. Hence (3) implies (1).\\ 
That (1) implies (3) follows from the definition of cohomology with local coefficients.
\end{proof}

The following theorem had already been expected by Fabien Morel, see \cite{ensprofin} \S 1.3.   
\begin{theorem}\label{modelstructure}
There is a left proper fibrantly generated model structure on $\hSh$ with the weak equivalences of Definition \ref{defnwe} for which $P$ is the set of generating fibrations and $Q$ is the set of generating trivial fibrations. The cofibrations  are the levelwise monomorphisms. We denote the homotopy category by $\hHh$.
\end{theorem}
\begin{proof}
In order to prove that there is a fibrantly generated model structure, we check the four conditions of the dual of Kan's Theorem 11.3.1 in \cite{hirsch}. Since we use cosmall instead of small objects, it suffices that $\hSh$ is closed under small limits and finite colimits. It is clear that the weak equivalences satisfy the 2-out-of-3 property and are closed under retracts. We denote by $Q$-cocell the subcategory of relative $Q$-cocell complexes consisting of limits of pullbacks of elements of $Q$. We write $P$-proj for the maps having the left lifting property with respect to all maps in $P$ and $P$-fib for the maps having the right lifting property with respect to all maps in $P$-proj. Now we check the remaining hypotheses of Kan's Theorem 11.3.1 in \cite{hirsch}.\\ 
(1) We have to show that the codomains of the maps in $P$ and $Q$ are cosmall relative to $P$-cocell and $Q$-cocell, respectively. This is clear for the terminal object $\ast$. It remains to check that the objects $K(M,n)$ and $BG$ are cosmall relative to $Q$-cocell. By definition of cosmallness we have to show that the canonical map 
$$f:\colim_{\alpha} \Hom_{\hSh}(Y_{\alpha},K(M,n)) \to \Hom_{\hSh}(\lim_{\alpha}Y_{\alpha},K(M,n))$$ is an isomorphism for some cardinal $\kappa$, where $Y_{\alpha}$ is any projective system whose indexing category is of cardinality $\kappa$ (and similarly for $BG$ instead of $K(M,n)$). By the definition of the spaces $K(M,n)$ (resp. $BG$) this map is equal to the map $\colim_{\alpha} Z^n(Y_{\alpha},M) \to Z^n(\lim_{\alpha}Y_{\alpha},M)$ (resp. with $M=G$ and $n=1$). But this map is already an isomorphism on the level of complexes $C^n$. For, Christensen and Isaksen have shown in \cite{prospectra}, Lemma 3.4, that  $\colim_{\alpha} \Hom_{\hEh}(Y_{\alpha}, M) \cong \Hom_{\hEh}(\lim_{\alpha}Y_{\alpha},M)$ (resp. $M=G$, $n=1$), since $\hEh$ is equivalent to the pro-category of finite sets.\\ 
(2) We have to show that every $Q$-fibration is both a $P$-fibration and a weak equivalence. Let $i:A\to B$ be a map in $P$-proj. As in the proof of Lemma \ref{lemme2}, this implies on the one hand that $i^{\ast}:Z^n(B,M)\to Z^n(A,M)$ (or $G$ instead of $M$ and $n=1$) is surjective and on the other hand that  $C^n(B,M)\to C^n(A,M)\times_{Z^{n+1}(A,M)}Z^{n+1}(B,M)$ is surjective for all $n\geq 0$ and every abelian finite group $M$ (or $M=G$ and $n=0$). It is easy to see that this implies that $i^{\ast}:C^n(B,M)\to C^n(A,M)$ (or $M=G$ and $n=1$) is surjective as well. Hence $i$ is an element in $Q$-proj. So $P$-proj $\subset$ $Q$-proj and hence $Q$-fib $\subset$ $P$-fib.\\ 
Furthermore, if $i$ is a monomorphism in each dimension, then $i$ is an element of $Q$-proj by Lemma \ref{lemme2}. Hence $Q$-fib is contained in the class of maps that have the right lifting property with respect to all monomorphisms. By Lemme 3 of \cite{ensprofin} this implies that the maps in $Q$-fib are simplicial homotopy equivalences and hence also weak equivalences.\\
(3) We have seen that every $P$-projective map is a $Q$-projective map. It remains to show $P$-proj $\subseteq W$. Let $f:A\to B$ be a map in $P$-proj. By the definition of $P$-proj, the maps $f^{\ast}:Z^1(B,G)\to Z^1(A,G)$ and $C^0(B,G)\to C^0(A,G)\times_{Z^{1}(A,G)}Z^{1}(B,G)$ are surjective. The latter implies that, if $f^{\ast}(\beta)$ is a boundary for an element $\beta \in Z^1(B,G)$, then $\beta$ is already a boundary. Hence $f$ induces an isomorphism $f^1:H^1(B;G) \cong H^1(A;G)$. Moreover, if $f$ has the left lifting property with respect to maps in $P$, then the same holds for its covering map $\tilde{f}$. The same argument for an abelian finite group $M$ instead of $G$ shows that for every $n\geq 1$ the map $\tilde{f}^n:H^n(\tilde{B};M) \to H^n(\tilde{A};M)$ is an isomorphism. For $n=0$, the left lifting property with respect to $K(S,0) \to \ast$ for any finite set $S$, not only shows that $H^0(f;S)=Z^0(f;S)$ is surjective but that it is also an isomorphism since any two liftings $B\to K(S,0)$ are simplicially homotopic and hence they agree. By Proposition \ref{whitehead}, this implies that $f$ is a weak equivalence.\\
(4) The remaining point is to show that $W\cap Q$-proj $\subseteq P$-proj. But this follows in almost the same way as we proved Lemma \ref{lemme2}.\\ 
This proves that we have found a fibrantly generated model structure on $\hSh$. It remains to show that it is left proper. In fact, the cofibrations are the maps in $Q$-proj and by Lemma \ref{lemme2} this class includes the monomorphisms. Hence every object in $\hSh$ is cofibrant which implies that the model structure is left proper. That the cofibrations are exactly the monomorphisms will be proved in the following lemma.
\end{proof}

\begin{lemma}\label{cofs}
A map in $\hSh$ is a cofibration if and only if it is a levelwise monomorphism. In particular, the maps $EG \to \ast$ and $L(M,n)\to \ast$ are trivial fibrations in $\hSh$ for every profinite group $G$, every abelian profinite group $M$ and every $n\geq 0$.
\end{lemma}
\begin{proof}
We have proven in the theorem that the class of cofibrations equals the class of maps having the left lifting property with respect to all maps in $Q$. We have seen that if $i$ is a monomorphism in each dimension, then $i$ is a cofibration by Lemma \ref{lemme2}.\\ 
So let $i:A \to B$ be a cofibration and suppose we have a map $f:A_n \to M$ in $\hEh$ for an abelian profinite group. Then the induced map $\Hom_{\hEh}(B_n,M/U) \to \Hom_{\hEh}(A_n,M/U)$ is surjective for every $n\geq 0$ and every open (and closed) normal subgroup $U$ of $M$, since $L(M/U,n)$ is an element of $Q$. Hence the set $N$ of pairs $(S,s)$ of closed subgroups $S$ of $M$ such that there is a lift $s:B_n\to M/S$ making the diagram
$$\xymatrix{
A_n \ar[d]_i  \ar[r]^{f/S} & M/S \\
B_n \ar[ur]_s}$$ 
commute contains the open normal subgroups of $M$. Then the same argument as in Lemma \ref{injectivity} shows that there is in fact a lift $B_n \to M$. Thus the induced map $i^{\ast}:\Hom_{\hEh}(B_n,M) \to \Hom_{\hEh}(A_n,M)$ is surjective for every abelian profinite group $M$ and every $n\geq 0$. The same argument works for an arbitrary profinite group $G$ and $n=0$.\\ 
By choosing $M$ to be the free profinite group on the set $A_n$, see e.g. \cite{ribes}, we see that $i$ must be injective in every dimension. Hence $Q$-proj equals the class of dimensionwise monomorphisms in $\hSh$.
\end{proof}

An important property of this model structure is the following homotopy invariance of limits in $\hSh$.
\begin{prop}\label{proplimits}
Let $f:X \to Y$ be a map of cofiltering diagrams $I \to \hSh$ of profinite spaces such that each $f_i:X_i \to Y_i$ is a weak equivalence. Then $\lim f_i:\lim X_i \to \lim Y_i$ is a weak equivalence in $\hSh$. In particular, the limit functor induces a functor on the homotopy category of $\hSh$.
\end{prop}
\begin{proof}
It follows from the fact that $\pi_0$ and $\pi_1$ commute with cofiltering limits that the profinite groupoids $\Pi(\lim_i X_i)$ and $\Pi(\lim_i Y_i)$ are equivalent. This also follows from the fact that cohomology with finite coefficients transforms limits to colimits, i.e. $H^0(\lim_i X_i;S)\cong \colim_i H^0(X_i;S)$ for a finite set $S$ and $H^1(\lim_i X_i;G)\cong \colim_i H^1(X_i;G)$ for a finite group $G$; hence the maps $H^0(f;S)$ and $H^1(f;G)$ are isomorphisms if each $H^0(f_i;S)$ and $H^1(f_i;G)$ is an isomorphism.
The cohomology with finite abelian local coefficient systems commutes with limits in the same way. By Proposition \ref{whitehead}, this implies that $\lim_i f_i$ is a weak equivalence.
\end{proof}

There is an obvious simplicial structure on $\hSh$, cf. \cite{ensprofin} and \cite{dehon}. The function complex $\hom_{\hSh}(X,Y)$ for $X,Y \in \hSh$ is the simplicial set defined in degree $n$ by $\Hom_{\hSh}(X\times \Delta[n],Y)$. It is characterized by the isomorphism
$$\Hom_{\hSh}(X\times W,Y) \cong \Hom_{\Sh}(W,\hom_{\hSh}(X,Y))$$
which is natural in the simplicial finite set $W$ and in $X,Y \in \hSh$. Moreover, if $j:A \to B$ is a cofibration and $q:X \to Y$ a fibration in $\hSh$, then the map 
$$\hom_{\hSh}(B,X) \longrightarrow \hom_{\hSh}(A,X)\times_{\hom_{\hSh}(A,Y)}\hom_{\hSh}(B,Y)$$
is a fibration of simplicial sets which is also a weak equivalence if $j$ or $q$ is one. In fact, by adjunction this statement is equivalent to that for every cofibration $i:V\to W$ in $\hSh$ the map
$$(A\times W)\cup_{A\times V} (B\times V) \longrightarrow B \times W$$
is a cofibration in $\hSh$ which is a weak equivalence if $j$ or $i$ is one in $\hSh$. The first point is clear and the second point follows from $\pi_1$ commuting with products, a van Kampen type theorem for profinite $\pi_1$, and the Mayer Vietoris long exact sequence for cohomology. Hence $\hSh$ is a simplicial model category in the sense of \cite{homalg}. In particular, for $X, Y \in \hSh$, $\hom_{\hSh}(X,Y)$ is fibrant in $\Sh$ if $Y$ is fibrant. Of course, there is also a simplicial model structure on the category $\hShp$ of pointed profinite spaces.\\
Furthermore, if $W$ is a simplicial set and $X$ is a profinite space, then by \cite{ensprofin} the function complex $\hom(W,X)$ has the natural structure of a profinite space as the cofiltering limit  of the simplicial finite sets $\hom(W_{\alpha},X/Q)$ where $W_{\alpha}$ runs through the simplicial finite subsets of $W$.\\
As an example we consider the simplicial finite set $S^1$, defined as $\Delta[1]$ modulo its boundary. For a pointed profinite space $X$, we denote its smash product $S^1\wedge X$ with $S^1$ by $\Sigma X$ and by $\Omega X$ the profinite space $\hom_{\hShp}(S^1,X)$. For $X,Y \in \hShp$, there is a natural bijection $\hom_{\hShp}(\Sigma X,Y)=\hom_{\hShp}(X, \Omega Y)$. Proposition \ref{whitehead} motivates the following definition for profinite higher homotopy groups.
\begin{defn}\label{defhomotopygroups}
Let $X$ be a pointed profinite space and let $RX$ be a functorial fibrant replacement of $X$ in the above model structure on $\hShp$. Then we define the {\em $n$th profinite homotopy group of $X$} for $n \geq 2$ to be the profinite group
$$\pi_n(X):=\pi_0(\Omega^n(RX)).$$
\end{defn}

It remains to show that the initial definition of the profinite fundamental group fits well with the definition of the higher homotopy groups, i.e. $\pi_1(X)\cong \pi_0(\Omega(RX))$. 
\begin{lemma}
For every $X \in \hSh$, the canonical map $\colim_Q \Ch_f(X/Q) \to \Ch_f(X)$ is bijective, where the colimit is taken over all $Q\in \Rh(X)$.
\end{lemma}
\begin{proof}
We have to show that any finite covering $E\to X$ of $X$ is induced by a finite covering of $X/Q$ via the quotient map $X \to X/Q$ for some $Q\in \Rh(X)$. We may assume that $E\to X$ is a Galois covering with finite Galois group $G$. Now $EG \to BG$ is the universal covering of $BG$. Hence by Proposition \ref{principalfibrations}, it suffices to note that $X\to BG$ is isomorphic to some quotient map $X\to X/Q$.
\end{proof}

\begin{prop}\label{pi1}
For a pointed profinite space $X$, the previously defined fundamental group $\pi_1X$ and the group $\pi_0\Omega X$ agree as profinite groups.
\end{prop}
\begin{proof}
The functors $\pi_1$, $\pi_0$ and $\Omega$ commute with cofiltering limits of fibrant objects by construction. The composed map 
$$X \stackrel{\cong}{\longrightarrow} \lim_{Q \in \Rh(X)} X/Q \stackrel{\simeq}{\longrightarrow} \lim_{Q\in \Rh(X)} R(X/Q)$$ 
is a weak equivalence by the homotopy invariance of cofiltering limits in $\hSh$ of Proposition \ref{proplimits}. Hence we may assume that $X$ is a fibrant simplicial finite set. In this case, $\pi_0\Omega X$ agrees with the usual $\pi_1 |X|$ of the underlying simplicial finite set of $X$ and hence we know that $\pi_0\Omega X$ is equal to the group of automorphisms of the universal covering of $X$. 
\end{proof}

Since $\hSh$ is a simplicial model category, for any profinite abelian group $M$, every $n$ and $X \in \hShp$ there is an isomorphism 
$$H^{n-q}(X;M)=\pi_q\hom_{\hShp}(X,K(M,n)),$$
where $\pi_q$ denotes the usual homotopy group of the simplicial set $\hom_{\hShp}(X,K(M,n))$. For an arbitrary profinite group $G$ there is a bijection of pointed sets 
$$H^{1}(X;G)=\pi_0\hom_{\hShp}(X,BG).$$
Let $X: I \to \Sh$ be a functor from a small cofiltering category $I$ to simplicial sets. By \cite{bouskan} XI \S 7.1, if each $X_i$ fibrant there is a spectral sequence involving derived limits
\begin{equation}\label{bouskanss}
E_2^{s,t}=\left \{ \begin{array}{l@{\quad \quad}l}
                     \lim_I^s\pi_t X_i & \mathrm{for}~0\leq s \leq t\\
                     0 & \mathrm{else}
                     \end{array} \right.
\end{equation}  
converging to $\pi_{s+t}\holim \, X$. Using \cite{bouskan} XI \S\S 4-7, one can construct a homotopy limit $\holim X \in \hSh$ for a small cofiltering category $I$ and a functor $X: I \to \hSh$. 
\begin{lemma}\label{holimlemma}
Let $X: I \to \hSh$ be a small cofiltering diagram such that each $X_i$ is fibrant in $\hSh$. Then there is a natural isomorphism
$$\pi_q(\holim X) \cong \lim_i \pi_q(X_i).$$
\end{lemma}
\begin{proof}
Since the underlying space of a fibrant profinite space is still fibrant in $\Sh$, the above spectral sequence exists also for a diagram of profinite spaces. But in this case all homotopy groups are profinite groups. Since the inverse limit functor is exact in the category of profinite groups, cf.~\cite{ribes} Proposition 2.2.4, the spectral sequence (\ref{bouskanss}) degenerates to a single row and implies the desired isomorphism.
\end{proof}

\begin{cor}\label{corlimholim}
Let $X:I \to \hSh$ be a cofiltering diagram  of profinite spaces such that each $X_i$ is fibrant. Then the natural map $\lim X \to \holim\, X$ is a weak equivalence in $\hSh$. 
\end{cor}
\begin{proof}
This follows directly from Lemma \ref{holimlemma} and the fact that the profinite homotopy groups commute with cofiltering limits of fibrant profinite spaces.
\end{proof}

\begin{cor}\label{propEM}
Let $I$ be a small cofiltering category and $\underline{M}$ (resp.~$\underline{G}$) be a functor from $I$ to the category of profinite abelian groups (resp.~profinite groups). Then the canonical map 
$K(\lim  \underline{M},n) \to \holim K(\underline{M},n)$ is a weak equivalence in $\hSh$ for every $n\geq 0$ (resp.~for $n=0,1$).  
\end{cor}
\begin{proof}
The lemma shows that $\pi_q \holim \,K(\underline{M},n)$ is equal to $M$ for $q=n$ and vanishes otherwise. Hence the canonical map $\lim K(\underline{M},n) \to \holim \,K(\underline{M},n)$ is a weak equivalence. The same argument holds for $G$ and $n=1$ using the construction of $\lim^1$ of \cite{bouskan} XI \S 6.   
\end{proof}

Since $\hom_{\hShp}(X,-)$ commutes with homotopy limits of fibrant objects, this result and the Bousfield-Kan spectral sequence (\ref{bouskanss}) imply the following result on the cohomology with profinite coefficients.
\begin{prop}\label{propfinite}
Let $f:X\to Y$ be a map in $\hSh$, let $G$ and $M$ be profinite groups and let $M$ be abelian. If  $H^1(f;G/U)$, resp.~$H^{\ast}(f;M/V)$, are isomorphisms for every normal subgroup of finite index $U \leq G$, resp.~$V \leq M$, then $H^1(f;G)$, resp.~$H^{\ast}(f;M)$, is an isomorphism. 
\end{prop} 
\begin{prop}
Let $M$ be a profinite group and suppose that the topology of $M$ has a basis consisting of a countable chain of open subgroups $M=U_0\geq U_1 \geq \ldots$. Then there is a natural short exact sequence for every $Y \in \hSh$ and every $i\geq 1$ ($i=1$ if $M$ is not abelian)
$$0 \to {\lim_n}^1 H^{i-1}(Y;M/U_n) \to H^i(Y;M) \to \lim_n H^i(Y;M/U_n) \to 0.$$  
\end{prop} 
\begin{proof}
This is the short exact sequence of \cite{bouskan} XI \S 7.4 for $X=K(M,i)$. Proposition \ref{propEM} identifies $\pi_i \holim \, X$ with $H^i(X;M)$.
\end{proof}

The hypothesis of the previous proposition is satisfied for example if $M$ is a finitely generated profinite abelian group. Finally, for a profinite group $G$ and a profinite abelian $G$-module $M=\lim_U M/U$, the isomorphism $$H^{n-q}(G;M)=\pi_q\hom_{\hSh/BG}(BG,K(M,n))$$ for the continuous cohomology of $G$ yields via (\ref{bouskanss}) and Proposition \ref{propEM} a natural spectral sequence for continuous group cohomology, $U$ running through the open normal subgroups of $M$:
$$E_2^{p,q}={\lim_U}^pH^q(G;M/U)\Rightarrow H^{p+q}(G;M).$$
Finally, we can generalize Lemma \ref{lemme2}.
\begin{prop}\label{EGtoBG}
For every profinite group $G$, every profinite abelian group $M$ and every $n\geq 0$, the canonical maps $EG \to BG$ and $L(M,n) \to K(M,n+1)$ are fibrations in $\hSh$.
\end{prop}
\begin{proof}
Let $X\hookrightarrow Y$ be a trivial cofibration in $\hSh$. Let $G$ be a profinite group, which is supposed to be abelian if $n\geq2$. We know by Lemma \ref{cofs} that the morphism of complexes $C^{\ast}(Y;G)\to C^{\ast}(X;G)$ is surjective. By assumption, it induces an isomorphism on the cohomology for every finite quotient $G/U$. Hence by Proposition \ref{propfinite} it also induces an isomorphism on cohomology with coefficients equal to $G$. Hence the maps $Z^n(Y;G) \to Z^n(X;G)$ and
$C^n(Y;G) \to C^n(X;G)\times_{Z^{n+1}(X;G)}Z^{n+1}(Y;G)$ are surjective and the maps $L(G,n) \to K(G,n+1)$ and $K(G,n) \to \ast$ have the desired right lifting property as in Lemma \ref{lemme2}.
\end{proof}

In terms of the homotopy category we can reformulate Proposition \ref{principalfibrations} as follows.
\begin{prop}
Let $G$ be a simplicial profinite group. For any profinite space $X$, the map
$$\theta:\Hom_{\hHh}(X,BG) \to \Phi^G(X),$$
sending the image of $f$ in $\Hom_{\hHh}(X,BG)$ to the pullback of $EG \to BG$ along $f$, is a bijection.
\end{prop}
\begin{cor}
Let $G$ be a simplicial profinite group and let $f:E\to X$ be a principal $G$-fibration. Then $f$ is also a fibration in $\hSh$.
\end{cor}
\begin{proof}
Since $X$ is cofibrant and $BG$ fibrant, the map $X\to BG$ in $\hHh$ corresponding to the principal fibration $E \to X$ is represented by a map in $\hSh$. Since $E \to X$ is the pullback of $EG\to BG$ under this map and since fibrations are stable under pullbacks, the assertion follows.
\end{proof}
The construction of the Serre spectral sequence of Dress in \cite{dress} can be easily translated to our profinite setting, see also \cite{dehon} \S 1.5. 
\begin{prop}\label{serress}
Let $B$ be a simply connected profinite space, let $f:E\to B$ be a fibration in $\hSh$ with fibre $F$ and let $M$ be an abelian profinite  group. Then there is a strongly convergent Serre spectral sequence 
$$E_2^{p,q}=H^p(B;H^q(F;M))\Rightarrow H^{p+q}(E;M).$$
\end{prop}
\begin{prop}\label{borel}
Let $G$ be a profinite group and let $p:E\to X$ be a principal $G$-fibration in $\hSh$. Then the canonical map $f:E\times_G EG \to X$ is a weak equivalence.
\end{prop}
\begin{proof}
Since $E\to X$ is locally trivial, see \cite{gj} V Lemma 2.5, it is also a covering of $X$ with fibre $G$. Hence we may assume that $p$ is a Galois covering of $X$ with $G=\mathrm{Aut}_X(E)$. It follows from the classification of coverings that there is a short exact sequence of groups 
$$1 \to \pi_1(E) \to \pi_1(X) \to G \to 1.$$ 
On the other hand, we have a canonical $G$-equivariant map $EG \times_{G}(E\times_G EG) = EG\times E \to E$ which induces an isomorphism on fundamental groups since $\pi_1(EG)$ is trivial. Now $EG \times_{G}(E\times_G EG)$ is a principal $G$-fibration on $E\times_G EG$ and we have a corresponding short exact sequence 
$$1 \to \pi_1(EG \times E) \to \pi_1(E\times_G EG) \to G \to 1.$$ 
Since the sequences are functorial, we conclude that $E\times_G EG \to X$ induces an isomorphism on fundamental groups. The remaining point to check follows from the Serre spectral sequence of Proposition \ref{serress} associated to the map of universal coverings of the fibration $E\times_G EG \to X$. 
\end{proof}

\subsection{Profinite completion of simplicial sets} 
We consider the category $\Sh$ of simplicial sets with the usual model structure of \cite{homalg}. We denote its homotopy category by $\Hh$.
\begin{prop}\label{adjcompletion}
1. The completion functor $\compl: \Sh \to \hSh$ preserves weak equivalences and cofibrations.\\
2. The forgetful functor $|\cdot|:\hSh \to \Sh$ preserves fibrations and weak equivalences between fibrant objects.\\
3. The induced completion functor $\compl: \Hh \to \hHh$ and the right derived functor $R|\cdot|:\hHh \to \Hh$ form a pair of adjoint functors.
\end{prop}
\begin{proof}
Let $f: X\to Y$ be a map of simplicial sets. If $f$ is a monomorphism and $x$ and $x'$ are two distinct $n$-simplices of $X/Q$ for some $Q\in \Rh(X)$, then there is a finite quotient $Y/R$, $R \in \Rh(Y)$ such that $f(x)$ and $f(x')$ are not equal. Hence $\hat{f}:\hat{X}\to \hat{Y}$ is a monomorphism.\\
If $f$ is a weak equivalence in $\Sh$, then $\pi_0(\hat{f})=\widehat{\pi_0(f)}$ and $\pi_1(\hat{f})=\widehat{\pi_1(f)}$ are isomorphisms for every basepoint by Proposition \ref{profinitecompletion}. Moreover, $H^n(f;\Mh)$ are isomorphisms for every finite abelian coefficient system $\Mh$ on $Y$ and every $n\geq 0$ by \cite{homalg}. Since the profinite completion of the universal covering of a space $X$ equals the universal profinite covering of the completion $\hat{X}$, we see that for a finite local system $\Mh$ the cohomologies $H^n(X;\Mh)$ and $H^n(\hat{X};\Mh)$ agree. Hence $\hat{f}$ is a weak equivalence in $\hSh$. The second and third assertion now follow from the first one since $\compl$ and $|\cdot|$ form a pair of adjoint functors.
\end{proof}

For a simplicial set $X$, we have seen that the continuous cohomology of $\hat{X}$ agrees with the cohomology of $X$ when the coefficients are finite. But for homotopical aspects, one should consider a fibrant replacement of $\hat{X}$ in $\hSh$ and could call this the profinite completion of $X$.

\subsection{Homology and the Hurewicz map}
We define the homology $H_{\ast}(X):=H_{\ast}(X;\hZ)$ of a profinite space $X$ to be the homology of the complex $C_{\ast}(X)$ consisting in degree $n$ of the profinite groups $C_n(X):=\hFab (X_n)$, the free abelian profinite group on the profinite set $X_n$. The differentials $d$ are the alternating sums $\sum_{i=0}^{n}d_i$ of the face maps $d_i$ of $X$. If $M$ is a profinite abelian group, then $H_{\ast}(X;M)$ is defined to be the homology of the complex $C_{\ast}(X;M):=C_{\ast}(X)\hat{\otimes}M$, where $\hat{\otimes}$ denotes the completed tensor product, see e.g.~\cite{ribes} \S 5.5. As for simplicial sets, there is a natural isomorphism of profinite groups, $H_0(X;\hat{\Z})=\hFab (\pi_0(X))$. For a pointed space $(X,\ast)$, we denote by $\tilde{H}_n(X;M)$ the reduced homology given by the complex $C_{\ast}(X;M)/C_{\ast}(\ast;M)$.\\
For a fibrant pointed profinite space $X$, by Proposition \ref{pi1} and by definition, the homotopy groups $\pi_n(X)$ are equal to the set $\Hom_{\hShp}(S^n,X)/\sim$ of maps modulo simplicial homotopy, where $S^n$ denotes the simplicial finite quotient $\Delta^n/\partial \Delta^n$. Hence an element $\alpha \in \pi_n(X)$ can be represented by an element $x \in X_n$. But we can view $x$ also as a cycle of $\tilde{C}_n(X)$ with homology class $[x]\in \tilde{H}_n(X;\hZ)$. One can show as in the case of simplicial sets that this correspondence $\alpha \mapsto [x]$ is well defined and is even a homomorphism of groups, cf.~\cite{may} \S 13. We call this map $h_n:\pi_n(X) \to \tilde{H}_n(X;\hZ)$ the Hurewicz map. 

\begin{prop}
Let $X$ be a connected pointed fibrant profinite space. The induced map $\overline{h}_1:\pi_1(X)/[\pi_1(X),\pi_1(X)] \to \tilde{H}_1(X;\hZ)$ is an isomorphism.
\end{prop}
\begin{proof}
This follows as in \cite{may} \S 13 or \cite{gj} III \S 3. Since $X$ has a strong deformation retract $Z$ that is reduced, i.e. $Z_0$ consists of a single element, we may assume that $X$ is reduced. Then $\tilde{C}_1(X)=\tilde{Z}_1(X)=\hFab (X_1-\{\ast\})$. The quotient $X_1\to X_1/\sim$ induces a natural epimorphism $j:\tilde{Z}_1(X) \to \pi_1/[\pi_1(X),\pi_1(X)]$. If $x\in \tilde{C}_2(X)$, then by the definition of $d$ and the definition of the group structure on $\pi_1(X)=\Hom_{\hShp}(S^1,X)/\sim$, we get $j\circ d(x)=0$. Thus $j$ induces a map $\overline{j}:\tilde{H}_1(X)\to \pi_1/[\pi_1(X),\pi_1(X)]$ and one checks easily that $\overline{j}$ and $\overline{h}$ are mutually inverse to each other.
\end{proof}

In exactly the same way as for simplicial sets, one proves the following Hurewicz theorem, see e.g. \cite{may} \S 13.
\begin{theorem}
Let $n\geq 1$ be an integer and let $X$ be a fibrant pointed profinite space with $\pi_q(X)=0$ for all $q<n$. Then the Hurewicz map $h:\pi_n(X)\to \tilde{H}_n(X;\hZ)$ is an isomorphism of profinite groups.
\end{theorem}

Let $X$ be a pointed space. By the universal property of profinite completion, there is canonical map $\widehat{\pi_n(X)}\to \pi_n(\hat{X})$ of profinite groups. We have seen in Proposition \ref{profinitecompletion} that this map is always an isomorphism for $n=1$. 
\begin{prop}\label{hatpin}
Let $X$ be a pointed simplicial set. Suppose that $\pi_q(X)=0$ for $q<n$. Then $\pi_n(\hat{X})$ is the profinite completion of $\pi_n(X)$, i.e.
$$\widehat{\pi_n(X)}\cong \pi_n(\hat{X}).$$
\end{prop}
\begin{proof}
This follows immediately from the Hurewicz theorem for profinite spaces.
\end{proof}

\subsection{Pro-$p$-model structures}
Morel has shown in \cite{ensprofin} that there is a $\Z/p$-model structure on $\hSh$ for every prime number $p$. The weak equivalences are maps that induce isomorphisms in the continuous cohomology with $\Zp$-coefficients. The cofibrations are the levelwise monomorphisms. It is also a fibrantly generated model structure. The minimal sets of generating fibrations and generating trivial fibrations are given by the canonical maps $L(\Z/p,n) \to K(\Z/p,n+1)$, $K(\Z/p,n)\to \ast$ and by the maps $L(\Z/p,n)\to \ast$ for every $n\geq 0$, respectively. Moreover, Morel proved that a principal $G$-fibration $E\to X$ is a fibration in this model structure for any pro-$p$-group $G$. In particular, the maps $EG\to BG$, $BG \to \ast$ are fibrations and the maps $EG\to \ast$ are trivial fibrations for any (nonabelian) pro-$p$-group $G$. Hence we would not get a different structure if we added for example the maps $EG \to BG$ for a nonabelian $p$-group to the generating sets of fibrations.\\
However, for the model structure of Theorem \ref{modelstructure} we cannot skip any map in the generating sets $P$ and $Q$. The arguments of \cite{ensprofin} rely on the fact that $\Zp$ is a field and that every pro-$p$-group has a $p$-central descending filtration by normal subgroups such that all subquotients are cyclic of order $p$. For a general profinite group there is not such a nice description.\\
One can describe the structure of \cite{ensprofin} as the left Bousfield localization of the model structure of Theorem \ref{modelstructure} with respect to the set of fibrant objects $K(\Z/p,n)$ for every $n\geq 0$. The homotopy groups $\pi_n(X):=\pi_0(\Omega^nX)$ of a fibrant profinite space $X$ for this model structure are pro-$p$-groups. For a simplicial set $X$, $\pi_1(\hat{X})$ is the pro-$p$-completion of $\pi_1(X)$. If $\pi_1(X)$ is finitely generated abelian, then $\pi_1(\hat{X})$ is isomorphic to $\Z_p\otimes_{\Z}\pi_1(X)$. 
\begin{remark}
The procedure for the proof of Theorem \ref{modelstructure} may be applied to every complete class $\Ch$ of finite groups. It suffices to replace the word {\em finite} by the appropriate additional property $\Ch$. For the class of $p$-groups the situation simplifies in the way indicated above.
\end{remark}
\subsection{Profinite completion of pro-spaces}
We define a completion functor $\compl:\pro-\Sh \to \hSh$ as the composite of two functors. First we apply $\compl :\Sh \to \hSh$ levelwise, then we take the limit in $\hSh$ of the underlying diagram. This completion functor unfortunately has no right adjoint. But it obviously preserves monomorphisms and weak equivalences by Propositions \ref{proplimits} and \ref{adjcompletion} when we equip pro-$\Sh$ with the model structure of \cite{modstruc}.\\ 
We compare the above constructions with Artin-Mazur's point of view. We denote by $\Hhf$ the full subcategory of $\Hh$ of spaces whose homotopy groups are all finite. If $X$ is either a space or a pro-space, then the functor $\Hhf \to \Eh$, $Y\mapsto \Hom_{\Hh}(Y,X)$ is pro-representable by \cite{artinmazur} Theorem 3.4. The corresponding pro-object $\hat{X}^{\mathrm{AM}}$ of $\Hhf$ is the Artin-Mazur profinite completion of $X$.\\
Now if $X$ is a profinite space we can naturally associate to it a pro-object $X^{\mathrm{AM}}$ in $\Hhf$. It is defined to be the functor $\Rh(RX) \to \Hhf$, $Q \mapsto RX/Q$, where $RX$ is a fibrant replacement of $X$ in $\hSh$. Since $RX/Q$ is finite in each degree, the same is true for $\Omega^k(RX/Q)$. Hence $\pi_0(\Omega^k(RX/Q))$ is finite for every $k\geq 0$. This functor $X \mapsto X^{\mathrm{AM}}$ sends weak equivalences to isomorphisms and hence factors through $\hHh$.\\
The composed functor $\Hh \to \hHh \to \mathrm{pro}-\Hhf$ is isomorphic to the Artin-Mazur completion of a space. In order to show this, it suffices to check that for a space $X$ the fundamental groups $\pi_1(\hat{X}^{\mathrm{AM}})$ and $\pi_1(\hat{X})$ are isomorphic as profinite groups and that cohomology with finite local coefficients agree. The first point follows from Proposition \ref{profinitecompletion} and \cite{artinmazur}, Corollary 3.7; and the second from the fact that cohomology with finite local coefficients transforms limits to colimits. In particular, this implies the following result.
\begin{prop}\label{artinmazurgroups}
Let $X$ be connected pointed pro-space. Then the homotopy pro-groups of the Artin-Mazur profinite completion $X^{\mathrm{AM}}$ and the profinite homotopy groups of $\hat{X}\in \hSh$ agree as profinite groups, i.e. for every $n\geq 1$
$$\pi_n(\hat{X})\cong \pi_n(X^{\mathrm{AM}}).$$
\end{prop}
Nevertheless, the categories $\hHh$ and $\mathrm{pro}-\Hhf$ are not equivalent as the example of Morel in \cite{ensprofin}, p. 368, shows.\\
Finally, we show that the continuous cohomology with profinite coefficients of the completion of a pro-space is equal to the continuous cohomology of Dwyer and Friedlander \cite{dwyfried} Definition 2.8, which we denote by $H^n_{\mathrm{DF}}(X;\Mh)$. 
\begin{prop}\label{contcohom}
Let $\Gamma$ be a profinite groupoid and let $X \to B\Gamma$ be a pro-space over $B\Gamma$. Let $\Mh$ be a profinite coefficient system on $\Gamma$. For every $n\geq 0$ there is a natural isomorphism induced by completion
$$H^n_{\Gamma}(\hat{X};\Mh)\cong H^n_{\mathrm{DF}}(X;\Mh).$$
\end{prop}
\begin{proof}
We may assume that $\Gamma$ is connected since we can decompose the cohomology into the product of the cohomology of the connected components. We denote by $\pi$ the profinite group $\Hom_{\Gamma}(\gamma,\gamma)$ of an object $\gamma \in \Gamma$. The continuous cohomology of \cite{dwyfried} is then given by 
$$H^n_{\mathrm{DF}}(X;\Mh)=\pi_0 \holim_j \colim_s \hom_{\Sh/B\pi}(X_s, K(M_j,n))$$
where $\Mh=\{M_j\}$ is given as an inverse system of finite $\pi$-modules and $\hom_{\Sh/B\pi}$ denotes the mapping space of spaces over $B\pi$. Because of the universal property of profinite completion and since each $K(M_j,n)$ is a simplicial finite set, we get a canonical identification of mapping spaces
$$\colim_s \hom_{\Sh/B\pi}(X_s, K(M_j,n))=\hom_{\hSh/B\pi}(\hat{X}, K(M_j,n)).$$
Furthermore, from Proposition \ref{propEM} we deduce
$$\pi_0 \holim_j \hom_{\hSh/B\pi}(\hat{X}, K(M_j,n))\cong \pi_0 \hom_{\hSh/B\pi}(\hat{X}, K(M,n))$$
where $M=\lim_j M_j$ is the profinite $\pi$-module corresponding to $\Mh$. Finally, the twisted cohomology of $\hat{X}$ is represented by the fibration $K(M,n) \to B\pi$ in $\hSh$, i.e. there is a canonical and natural identification 
$$H_{\pi}^n(\hat{X};M)= \pi_0\hom_{\hSh/B\pi}(\hat{X},K(M,n)).$$ 
This series of isomorphisms now yields the proof of the assertion.
\end{proof}

\subsection{Stable profinite homotopy theory}
\begin{defn}\label{defspectrum}
A {\rm profinite spectrum} $X$ consists of a sequence $X_n \in \hShp$ of pointed profinite spaces for $n\geq0$ and maps $\sigma_n:S^1 \wedge X_n \to X_{n+1}$ in $\hShp$. A {\rm morphism} $f:X \to Y$ of spectra consists of maps $f_n:X_n \to Y_n$ in $\hShp$ for $n\geq0$ such that $\sigma_n(1\wedge f_n)=f_{n+1}\sigma_n$. We denote by $\hSp$ the corresponding category and call it the category of {\rm profinite spectra}. A spectrum $E\in \hSp$ is called an {\rm $\Omega$-spectrum} if each $E_n$  is fibrant and the adjoint structure maps $E_n \to \Omega E_{n+1}$ are weak equivalences for all $n\geq 0$.
\end{defn}
The suspension $\Sigma^{\infty}: \hShp \to \hSp$ sends a profinite space $X$ to the spectrum given in degree $n$ by $S^n\wedge X$. Starting with the model structure on $\hSh$ of Theorem \ref{modelstructure}, the localization theorem of \cite{etalecob} yields the following result.
\begin{theorem}\label{thmprostable}
There is a {\rm stable model structure on $\hSp$} for which the prolongation of the suspension functor is a Quillen equivalence. The corresponding stable homotopy category is denoted by $\hSHh$. In particular, the {\rm stable equivalences} are the maps that induce an isomorphism on all generalized cohomology theories, represented by profinite $\Omega$-spectra; the {\rm stable cofibrations} are the maps $i:X \to Y$ such that $i_0$ and the induced maps $j_n: X_n \amalg_{S^1\wedge X_{n-1}}S^1 \wedge Y_{n-1} \to Y_n$ are monomorphisms for all $n$; the {\rm stable fibrations} are the maps with the right lifting property with respect to all maps that are both stable equivalences and stable cofibrations.
\end{theorem}
Since the suspension is compatible with profinite completion, there is an analogous statement as in Proposition \ref{adjcompletion} on the pair of adjoint functors consisting of the levelwise completion functor $\Sp \to \hSp$ and the forgetful functor where $\Sp$ denotes the category of simplicial spectra with the Bousfield-Friedlander model structure \cite{bousfried}. For a profinite spectrum $X$, we define the stable homotopy groups $\pi^s_n(X)$ to be the set $\Hom_{\hSHh}(S^n,X)$ of maps in the stable homotopy category. Since $\Omega$ preserves fibrations in $\hShp$, Proposition \ref{hatpin} implies the following characterization of the stable homotopy groups of the completion of a spectrum.
\begin{prop}
Let $E$ be a connected simplicial spectrum, then the canonical map $\widehat{\pi^s_n(E)} \to \pi^s_n(\hat{E})$ is an isomorphism. In particular, if each $\pi^s_n(E)$ is finitely generated, then $\pi^s_n(\hat{E})\cong \hZ \otimes_{\Z} \pi^s_n(E)$.
\end{prop}
Examples for the last assertion are Eilenberg-MacLane spectra $HM$ for a finitely generated abelian group $M$ or the simplicial spectrum of complex cobordism $MU$. Finally, the composition of completion and homotopy limit defines also a functor from pro-spectra to profinite spectra that induces a functor from the homotopy category of pro-spectra of Christensen and Isaksen in \cite{prospectra} to $\hSHh$.

\section{Applications in \'etale homotopy theory}

\subsection{Profinite \'etale homotopy groups}

The construction of the \'etale topological type functor $\Et$ from locally noetherian schemes to pro-spaces is due to Artin-Mazur and Friedlander. We refer the reader to \cite{fried} and \cite{a1real} for a detailed discussion of the category of rigid hypercoverings and rigid pullbacks. Let $X$ be a locally noetherian scheme. The {\rm \'etale topological type of $X$} is defined to be the pro-simplicial set $\Et X := \mathrm{Re} \circ \pi :HRR(X) \to \Sh$ sending a rigid hypercovering $U$ of $X$ to the simplicial set of connected components of $U$. If $f: X \to Y$ is a map of locally noetherian schemes, then the strict map $\Et f:\Et X \to \Et Y$ is given by the functor $f^{\ast}:HRR(Y) \to HRR(X)$ and the natural transformation $\Et X \circ \Et f \to \Et Y$. In order to get a profinite space, we compose $\Et$ with the completion from pro-spaces to $\hSh$ and denote this functor by $\hEt$.\\
Let $X$ be a pointed connected locally noetherian scheme. For $k>0$, let $\pi_k(\Et X)$ be the pro-group, defined by the functor $\pi_k \circ \Et$ from $HRR(X)$ to the category of groups as in \cite{fried}. For $k=1$ and $G$ a group, Friedlander has shown in \cite{fried}, 5.6, that the set of isomorphism classes of principal $G$-fibrations over $X$ is isomorphic to the set $\Hom(\pi_1(\Et X),G)$ of homomorphisms of pro-groups. Furthermore, the locally constant \'etale sheaves on $X$, whose stalks are isomorphic to a set $S$, are in $1$-$1$ correspondence to local coefficient systems on $\Et X$ with fibers isomorphic to $S$. This shows that $\pi_1(\Et X)$ is equal to the enlarged fundamental pro-group of \cite{sga3}, Expos\'e X, \S 6. It agrees with the profinite \'etale fundamental group of \cite{sga1} if $X$ is geometrically unibranched (e.g. if $X$ is normal). More generally, Friedlander shows in \cite{fried}, Theorem 7.3, that $\pi_k(\Et X)$ is a profinite group for $k>0$ when $X$ is connected and geometrically unibranched.\\
These arguments are easily transfered to the profinite setting. For a topological group $G$ there is a bijection between the set of isomorphism classes of principal $G$-fibrations over $X$ and the set of isomorphism classes of principal $G$-fibrations over $\hEt X$; and hence there is a bijection between the set of isomorphism classes of locally constant \'etale sheaves of profinite sets on $X$ and the set of isomorphism classes of profinite local coefficient systems on $\hEt X$. Furthermore, we see that the finite coverings of $\hEt X$ are in $1$-$1$-correspondence with the finite \'etale coverings of $X$. This implies that the finite quotients of $\pi_1(\hEt X)$ correspond to the finite \'etale coverings of $X$. Hence the profinite group $\pi_1(\hEt X)$ agrees with the usual profinite fundamental group $\pi_1^{\et}(X)$ of \cite{sga1}.\\
If $F$ is a locally constant \'etale sheaf on $X$, we denote by $F$ also the corresponding local system on $\Et X$ respectively $\hEt X$. Friedlander has shown in \cite{fried} Proposition 5.9, that the groups $H_{\et}^{\ast}(X;F)$ and $H^{\ast}(\Et X;F)$ are equal. Moreover, by the definition of $\hEt X$, the cohomology groups $H^{\ast}(\hEt X;F)$ and $H^{\ast}(\Et X;F)$ coincide if $F$ is finite. But if $F$ is not finite, $H^{\ast}(\hEt X;F)$ does not in general agree with the usual \'etale cohomology groups $H_{\et}^{\ast}(X;F)$ any more. But since, for example, $\ell$-adic cohomology $H^{i}(X;\Z_{\ell}(j)):=\lim_n H_{\et}^{i}(X;\Z/\ell^n(j))$ of a scheme $X$ has good properties only if the \'etale cohomology groups $H_{\et}^{i}(X;\Z/\ell^n(j))$ are finite, this may not be considered as a problem. In fact, it turns out that, if $F$ is profinite, $H^{\ast}(\hEt X;F)$ is the continuous \'etale cohomology of $X$, which is a more sophisticated version of \'etale cohomology for inverse systems of  coefficient sheaves defined by Jannsen in \cite{jannsen}. 
\begin{prop}
Let $X$ be a locally noetherian scheme and let $F$ be a locally constant \'etale sheaf on $X$ whose stalks are profinite abelian groups. The cohomology $H^{\ast}(\hEt X;F)$ of $\hEt X$ with coefficients in the local system corresponding to $F$ coincides with the continuous \'etale cohomology $H^{\ast}_{\mathrm{cont}}(X;F)$ of Jannsen and of Dwyer and Friedlander. In particular, $H^{\ast}(\hEt X;\Z_{\ell}(j))$ equals $H_{\mathrm{cont}}^{\ast}(X;\Z_{\ell}(j))$ of \cite{jannsen}.
\end{prop}
\begin{proof}
This follows immediately from Lemma 3.30 of \cite{jannsen} and Proposition \ref{contcohom}, since $\Et X$ is a pro-space over the profinite groupoid $\Pi_1^{\et}(X)=\Pi(\hEt X)$. 
\end{proof}

The relation to the usual $\ell$-adic cohomology of a locally noetherian scheme $X$ is given by the exact sequence, as in \cite{jannsen},
$$0 \to {\lim_n}^1 H_{\et}^{i-1}(X;\Z/\ell^n(j)) \to H^i(\hEt X;\Z_{\ell}(j)) \to \lim_n H_{\et}^i(X;\Z/\ell^n(j)) \to 0.$$ 
If $p:Y \to X$ is a Galois covering with profinite Galois group $G$ and $F$ is a locally constant profinite \'etale sheaf on $X$, then there is a spectral sequence of continuous cohomology groups 
$$E_2^{p,q}=H^p(G;H^q(\hEt Y;p^{\ast}F)) \Rightarrow H^{p+q}(\hEt X;F).$$
The previous discussion shows that $\hEt X$ is a good rigid model for the profinite homotopy type of a scheme. Since the \'etale fundamental group $\pi_1^{\et}(X)$ is always a profinite group and is equal to $\pi_1(\hEt X)$, we make the following definition.
\begin{defn} 
For a locally noetherian simplicial scheme $X$, a geometric point $x$ and $n\geq 2$, we define the {\rm profinite \'etale homotopy groups} $\pi_n^{\et}(X,x)$ of $X$ to be the profinite groups $\pi_n(\hEt X,x)$.
\end{defn}
By Proposition \ref{artinmazurgroups}, these profinite homotopy groups agree with the profinitely completed \'etale homotopy groups of Artin and Mazur in \cite{artinmazur}.

\subsection{Example: $\hEt k$ and the absolute Galois group}
Let $k$ be a fixed field. For a Galois extension $L/k$ we denote by $G(L/k)$ its Galois group. For a separable closure $\ok$ of $k$ we write $G_k$ for the absolute Galois group $G(\ok/k)$ of $k$. We want to determine the well known homotopy type of $\hEt k$ as $BG_k$. We denote by $RC(k)$ the category of rigid coverings and by $HRR(k)$ the category of rigid hypercoverings of $\Spec k$. There are two canonical functors between them. On the one hand there is the restriction functor sending a hypercovering $U$ to $U_0$; on the other hand we can send a rigid covering $U \to \Spec k$ to the hypercovering $\cosk_0^k(U)$. When we consider the rigid cover $\Spec L \to \Spec k$ associated to a finite Galois extension $L/k$, the simplicial set of connected components of the corresponding hypercover $\pi_0(\cosk_0^k(L))$ is equal to the classifying space $BG(L/k)$. Hence when we consider the functor $\pi_0\circ \cosk_0^k(-):RC(k)\to \Sh$ from rigid coverings of $k$, we see that it takes values in simplicial profinite sets $\hSh$ and that the limit over all the rigid coverings equals $BG_k$.\\
Furthermore, the restriction functor defines a monomorphism $g:BG_k \hookrightarrow \hEt k$; and the coskeleton $\cosk_0^k$ defines a map in the other direction $f:\hEt k \to BG_k$. In particular, $f$ and $g$ are simplicial homotopy equivalences and $g$ is a cofibration. For, it is clear that $f\circ g$ equals the identity of $BG_k$, whereas the composite $g\circ f$ is homotopic to the identity of $\hEt k$ by \cite{fried} Prop. 8.2.\\ 
For a $k$-scheme $X$ we write $X_L:= X\otimes_k L$. The group $G:=G(L/k)$ acts via Galois-automorphisms $1 \otimes \sigma$ on $X_L$. The canonical map $X_L \to X$ is an \'etale Galois covering of $X$ with Galois group $G$. Hence $\hEt X_L \to \hEt X$ is a principal $G$-fibration in $\hSh$. This implies that the canonical map $\hEt X_L \times_{G} EG \to \hEt X$ is a weak equivalence by Proposition \ref{borel} and that the cohomology of $\hEt X$ is equal to the Borel-cohomology of $\hEt X_L$ under the action of $G(L/k)$.\\
An interesting and well known example is the \'etale realization of a finite field $\F_q$ with $q$ elements. Its Galois group is $\hZ$ and hence $\hEt \F_q=K(\hZ,1)$. In fact, there is a canonical homotopy equivalence $S^1 \to \hEt \F_q$ sending the generator of $\pi_1S^1=\hZ$ to the Frobenius map.\\ 
Finally, for a base field $k$, $\hEt \Pro_{\ok}^1$ is an Eilenberg MacLane $G_k$-space $K(\hZ(1),2)$ with the canonical $G_k$-action on the profinite group $\hZ(1)=\mu(\ok)$ of all roots of unity in $\ok$. More examples may be found in \cite{topmodels}.
\subsection{Etale realization of the flasque motivic model structure}
The construction of the $\A^1$-homotopy category of schemes gave rise to the question if this functor may be enlarged to the category of motivic spaces. This has been answered independently by Dugger-Isaksen and Schmidt. The latter one constructed a geometric functor to the category pro-$\Hh$ as Artin-Mazur. The idea of Dugger and Isaksen for the less intuitive extension of $\Et X$ is that it should be the usual $\Et X$ on a representable presheaf $X$ and should preserve colimits and the simplicial structure. Isaksen then showed that $\Et$ induces a left Quillen functor on the projective model structure on simplicial presheaves.\\ 
Let $S$ be a base scheme and let $\Sm/S$ be the category of smooth quasi-projective schemes of finite type over $S$. Isaksen has shown in \cite{flasque} that the flasque model structure on presheaves on $\Sm/S$ is a good model for the $\A^1$-homotopy category, in particular for the construction of the stable motivic homotopy category, since $\Pro^1$ is a flasque cofibrant space. It was shown in \cite{etalecob} that $\hEt$ yields a stable \'etale realization functor of the stable motivic homotopy category. In order to simplify the constructions, we will deduce from \cite{a1real} that $\hEt$ induces a derived functor on the flasque motivic structure as well. 
\begin{lemma}\label{immersions}
Let $i:Y \hookrightarrow X$ be an open or closed immersion of locally noetherian schemes. Then $\Et i:\Et Y \hookrightarrow \Et X$ is a monomorphism of pro-simplicial sets.
\end{lemma}
\begin{proof}
The map $\Et i$ is by definition the map of pro-simplicial sets which is induced by the natural transformation of indexing categories given by the rigid pullback functor $i^{\ast}:HRR(X) \to HRR(Y)$. If $\coprod_{x}U_x,x \to X$ is a rigid cover of $X$ indexed by the geometric points of $X$, the connected component $(i^{\ast}U)_y$ of the pullback cover of $Y$ at the geometric point $y$ of $Y$ is sent to the connected component $U_{i(y)}$ over $X$. If $i$ is an open, resp. closed, immersion, this map on the sets of connected components is obviously injective. Since open, resp. closed, immersions are stable under base change and since for a rigid hypercover $V$ of $X$, 
$(\cosk_{t-1}^{X_s}V)_t$ is a finite fiber product involving the $V_{s,r}$'s for $r\leq t$ and $X_s$, we deduce that the induced map $\Et i$ is a monomorphism as well. 
\end{proof}

\begin{theorem}
The \'etale topological type functor $\Et$ is a left Quillen functor on the Nisnevich (resp. \'etale) local flasque model structure on simplicial presheaves on $\Sm/S$ to the model structure on pro-simplicial sets of \cite{modstruc}. 
\end{theorem}
\begin{proof}
Let $i:\cup_{k=1}^nU_n \to X$ be an acceptable monomorphism as defined in \cite{flasque}. We have to show that $\Et i$ is a monomorphism in pro-$\Sh$. By compatibility with colimits, $\Et \cup_{k=1}^nU_n$ is the coequalizer of the diagram
$$\Et(\coprod_{k,k'}U_k\times_X U_{k'}) \to \Et(\coprod_k U_k).$$ Hence in order to show that $\Et i$ is a cofibration, it is enough to show that $\Et(U_k\times_X U_{k'}) \to \Et U_k$ is a monomorphism for every $k$ and $k'$. But since open immersions of schemes are stable under base change, this follows from Lemma \ref{immersions}. Moreover, $\Et$ is compatible with pushout products and hence $\Et$ sends the generating (trivial) cofibrations in the flasque model structure on presheaves to (trivial) cofibrations of pro-simplicial sets.\\
The Nisnevich (resp. \'etale) local flasque model structure is constructed via the left Bousfield localizations with respect to all Nisnevich (resp. \'etale) hypercovers. In fact, Isaksen shows that $\Et$ sends Nisnevich and \'etale hypercovers to weak equivalences in pro-$\Sh$, \cite{a1real}, Theorem 12.   
\end{proof}

\begin{cor}
If $\car S=\{0\}$, the functor $\hEt$ also induces a total left derived functor from the motivic flasque model structure on $\sPreS$ to $\hSh$ with the general model structure of Theorem \ref{modelstructure}.\\
If $\car S$ contains a prime $p>0$, it induces a total derived functor to $\hSh$ with the $\Zl$-model structure of \cite{ensprofin} for any prime number $\ell$ which is prime to $\car S$.
\end{cor}
\begin{proof}
Since $\Et$ induces a total left derived functor on the local structures and since the above completion induces a functor on the homotopy categories, it is clear that $\hEt$ induces the desired functor for the local structures.\\
It remains to check that $\hEt$ sends the projection $p_X:\A^1 \times_k X \to X$ to a weak equivalence in $\hSh$ for all $X\in \Sm/k$. If $\car k=0$, $p_X$ induces an isomorphism on \'etale cohomology $H_{\et}^{\ast}(-;F)$ for all torsion sheaves $F$, see e.g. \cite{milne} VI, Corollary 4.20, and an isomorphism on \'etale fundamental groups. Hence in this case $\hEt p_X$ is a weak equivalence in $\hSh$.\\
If $\car k=p >0$, $\pi_X$ induces an isomorphism on \'etale cohomology $H_{\et}^{\ast}(-;F)$ for all torsion sheaves $F$ whose torsion is prime to $p$. But it does not induce an isomorphism on \'etale fundamental groups. Hence in this case $\hEt p_X$ is only a weak equivalence in the mod $\ell$-model structure on $\hSh$. 
\end{proof}

The space $\Pro_k^1$ pointed at $\infty$ is flasque cofibrant and may be used as in \cite{flasque} to construct the stable motivic homotopy category of $\Pro^1$-spectra starting from the flasque model structure on motivic spaces.  As in \cite{etalecob} this yields a stable motivic realization functor without taking a cofibrant replacement of $\Pro^1$.
\begin{cor}
Let $k$ be a field with $\car k=0$. The functor $\hEt$ induces a stable \'etale realization functor of the stable motivic homotopy category to the homotopy category $\hSHh$ of Theorem \ref{thmprostable}.
\end{cor}
\begin{remark}
1. For a pointed presheaf $\Xh$ on $\Sm/S$ there are profinite \'etale homotopy groups for every $n\geq 1$: $\pi_n^{\et}(\Xh):=\pi_n(\hEt \Xh)$.\\ 
2. All statements of this subsection also hold for the closed (motivic) model structure of \cite{ppr}.
\end{remark}
\bibliographystyle{amsplain}

Gereon Quick, Mathematisches Institut, Universit\"at M\"unster, Einsteinstr. 62, D-48149 M\"unster\\
E-mail address: gquick@math.uni-muenster.de\\
Homepage: www.math.uni-muenster.de/u/gquick

\end{document}